\documentclass[a4paper,10pt]{article}

\AtBeginDocument{%
  }

\def\longtitle{Enabling~Research~through~the~SCIP~Optimization~Suite~8.0}
\def\shortfunding{The work for this article has been partly conducted within the
  \emph{Research Campus MODAL} funded by the German Federal Ministry of
  Education and Research (BMBF grant number 05M14ZAM) and has received funding from
  the European Union's Horizon 2020 research and innovation programme under
  grant agreement No 773897. It has also been
  partly supported by
  the German Research Foundation (DFG) within the
  Collaborative Research Center 805, Project A4,
  and the EXPRESS project of the priority program CoSIP (DFG-SPP 1798),
  the German Research Foundation (DFG) within the project HPO-NAVI (project number 391087700).
  }

\usepackage{xspace}
\usepackage{amsmath}
\usepackage{amssymb}
\usepackage{mathtools}
\usepackage{xcolor}
\usepackage{hyperref}
\usepackage{graphicx}
\usepackage{booktabs}
\usepackage[numbers]{natbib}

\makeatletter
\newcommand*\rel@kern[1]{\kern#1\dimexpr\macc@kerna}
\newcommand*\widebar[1]{%
  \begingroup
  \def\mathaccent##1##2{%
    \rel@kern{0.8}%
    \overline{\rel@kern{-0.8}\macc@nucleus\rel@kern{0.2}}%
    \rel@kern{-0.2}%
  }%
  \macc@depth\@ne
  \let\math@bgroup\@empty \let\math@egroup\macc@set@skewchar
  \mathsurround\z@ \frozen@everymath{\mathgroup\macc@group\relax}%
  \macc@set@skewchar\relax
  \let\mathaccentV\macc@nested@a
  \macc@nested@a\relax111{#1}%
  \endgroup
}
\makeatother

\newcommand{\LP}{{LP}\xspace}
\newcommand{\LPs}{{LPs}\xspace}
\newcommand{\CIP}{{CIP}\xspace}
\newcommand{\CIPs}{{CIPs}\xspace}

\newcommand{\MILP}{{MILP}\xspace}
\newcommand{\MILPs}{{MILPs}\xspace}

\newcommand{\MINLP}{{MINLP}\xspace}
\newcommand{\MINLPs}{{MINLPs}\xspace}

\newcommand{\MISDP}{{MISDP}\xspace}

\newcommand{\T}{\top}

\newcommand{\defi}{\coloneqq}

\newcommand{\linobj}{c}

\newcommand{\linmatrix}{A}
\newcommand{\nonlincons}{g}
\newcommand{\rhs}{b}

\newcommand{\low}[1]{\underline{#1}}
\newcommand{\upp}[1]{\overline{#1}}
\newcommand{\consindex}{\mathcal{M}}
\newcommand{\varindex}{\mathcal{N}}
\newcommand{\intvarindex}{\mathcal{I}}
\newcommand{\NP}{\ensuremath{\mathcal{NP}}\xspace}

\newcommand{\N}{\mathbb{N}}
\newcommand{\R}{\mathbb{R}}

\newcommand{\Z}{\mathbb{Z}}
\newcommand{\Rinf}{\ensuremath{\widebar{\mathbb{R}}}\xspace}

\usepackage{siunitx}
\newcommand{\cleaninst}{all}
\newcommand{\affected}{affected}
\newcommand{\alloptimal}{{both-solved}\xspace}
\newcommand{\difftimeouts}{{diff-timeouts}\xspace}
\newcommand{\bracket}[2]{[#1,#2]}

\definecolor{c1}{HTML}{000060}
\definecolor{c2}{HTML}{0000FF}
\definecolor{c3}{HTML}{36648B}
\definecolor{c4}{HTML}{4682B4}
\definecolor{c5}{HTML}{5CACEE}
\definecolor{c6}{HTML}{00FFFF}
\definecolor{c7}{HTML}{008888}
\definecolor{c8}{HTML}{00DD99}
\definecolor{c9}{HTML}{527B10}
\definecolor{c10}{HTML}{7BC618}
\definecolor{c11}{HTML}{33AA00}

\definecolor{scipoldcol}{HTML}{36648B}
\definecolor{scipnewcol}{HTML}{7BC618}

\newcommand{\solver}[1]{\textsc{#1}\xspace}
\newcommand{\scipopt}{\scip Optimization Suite\xspace}
\newcommand{\scipversion}{8.0}

\newcommand{\scipoptv}{\scipopt~\scipversion\xspace}
\newcommand{\scip}{\solver{SCIP}}
\newcommand{\scipv}{\solver{SCIP}~\scipversion\xspace}
\newcommand{\soplex}{\solver{SoPlex}}
\newcommand{\soplexversion}{6.0}
\newcommand{\soplexv}{\solver{SoPlex}~\soplexversion\xspace}
\newcommand{\papilo}{\solver{PaPILO}}
\newcommand{\papiloversion}{2.0}
\newcommand{\papilov}{\solver{PaPILO}~\papiloversion\xspace}
\newcommand{\zimpl}{\solver{Zimpl}} 

\newcommand{\ug}{\solver{UG}}
\newcommand{\presollib}{\solver{PaPILO}}

\newcommand{\gcg}{\solver{GCG}}
\newcommand{\gcgversion}{3.5}

\newcommand{\scipsdp}{\solver{SCIP-SDP}}

\newcommand{\scipjack}{\solver{SCIP-Jack}}

\newcommand{\scipjackversion}{2.0}

\newcommand{\cplex}{\solver{CPLEX}}
\newcommand{\ipopt}{\solver{Ipopt}}
\newcommand{\cppad}{\solver{CppAD}}

\newcommand{\gurobi}{\solver{Gurobi}}

\newcommand{\nbsc}[1]{\mbox{#1}\xspace}
\newcommand{\miplib}{\nbsc{MIPLIB}}

\newcommand{\coral}{\nbsc{COR@L}}

\newcommand{\minlplibtwo}{\nbsc{MINLPLib}}

\definecolor{darkgreen}{HTML}{008800}

\newcommand{\fa}{\text{ for all }}

\newcommand{\perm}{\gamma}
\newcommand{\inv}[1]{{#1}^{-1}}

\newcommand{\group}{\Gamma}

\newcommand{\bliss}{\solver{bliss}}

\newcommand{\myand}{$\cdot$\xspace}
\newcommand{\myorcidlink}[1]{\,\href{https://orcid.org/#1}{\raisebox{-0.45ex}{\includegraphics[width=1.8ex]{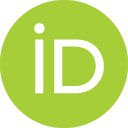}}}}

\begin{document}

\title{\longtitle}

\author{\small
  Ksenia Bestuzheva\protect\myorcidlink{0000-0002-7018-7099} \myand
  Mathieu Besançon\protect\myorcidlink{0000-0002-6284-3033} \myand
  Wei-Kun Chen \protect\myorcidlink{0000-0003-4147-1346} \and\small
  Antonia Chmiela\protect\myorcidlink{0000-0002-4809-2958} \myand
  Tim Donkiewicz\protect\myorcidlink{0000-0002-5721-3563} \myand
  Jasper van Doornmalen\protect\myorcidlink{0000-0002-2494-0705} \and\small
  Leon Eifler\protect\myorcidlink{0000-0003-0245-9344} \myand
  Oliver Gaul\protect\myorcidlink{0000-0002-2131-1911} \myand
  Gerald Gamrath\protect\myorcidlink{0000-0001-6141-5937} \myand
  Ambros Gleixner\protect\myorcidlink{0000-0003-0391-5903} \and\small
  Leona Gottwald\protect\myorcidlink{0000-0002-8894-5011} \myand
  Christoph Graczyk\protect\myorcidlink{0000-0001-8990-9912} \myand
  Katrin Halbig\protect\myorcidlink{0000-0002-8730-3447} \myand
  Alexander Hoen\protect\myorcidlink{0000-0003-1065-1651} \and\small
  Christopher Hojny\protect\myorcidlink{0000-0002-5324-8996} \myand
  Rolf van der Hulst\protect\myorcidlink{0000-0002-5941-3016} \myand
  Thorsten Koch\protect\myorcidlink{0000-0002-1967-0077} \and\small
  Marco L\"ubbecke\protect\myorcidlink{0000-0002-2635-0522} \myand
  Stephen J.~Maher\protect\myorcidlink{0000-0003-3773-6882} \myand
  Frederic Matter\protect\myorcidlink{0000-0002-0499-1820} \myand
  Erik M\"uhmer\protect\myorcidlink{0000-0003-1114-3800} \and\small
  Benjamin M\"uller\protect\myorcidlink{0000-0002-4463-2873} \myand
  Marc E. Pfetsch\protect\myorcidlink{0000-0002-0947-7193} \myand
  Daniel Rehfeldt\protect\myorcidlink{0000-0002-2877-074X} \myand
  Steffan Schlein\protect\myorcidlink{0009-0003-8322-7226} \and\small
  Franziska Schlösser\protect\myorcidlink{0000-0002-3716-5031} \myand
  Felipe Serrano\protect\myorcidlink{0000-0002-7892-3951} \myand
  Yuji Shinano\protect\myorcidlink{0000-0002-2902-882X} \myand
  Boro Sofranac\protect\myorcidlink{0000-0003-2252-9469} \and\small
  Mark Turner\protect\myorcidlink{0000-0001-7270-1496} \myand
  Stefan Vigerske\protect\myorcidlink{0009-0001-2262-0601} \myand
  Fabian Wegscheider\protect\myorcidlink{0009-0000-8100-6751} \myand
  Philipp Wellner\protect\myorcidlink{0009-0001-1109-3877} \and\small
  Dieter Weninger\protect\myorcidlink{0000-0002-1333-8591} \myand
  Jakob Witzig\protect\myorcidlink{0000-0003-2698-0767}%
  \thanks{Extended author information is available at the end of the paper.
    \shortfunding}}

\maketitle

\begin{abstract}
The \scipopt provides a collection of software packages for mathematical optimization centered
around the constraint integer programming framework \scip.
The focus of this paper is on the role of the \scipopt in supporting research.
\scip's main design principles are discussed, followed by a presentation of the latest performance
improvements and developments in version 8.0,
which serve both as examples of \scip's application as a research tool and as a platform for further
developments.
Further, the paper gives an overview of interfaces to other programming and modeling languages,
new features that expand the possibilities
for user interaction with the framework,
and the latest developments in several extensions built upon \scip.
\end{abstract}

\paragraph{\bf Keywords} Constraint integer programming
$\cdot$ linear programming
$\cdot$ mixed-integer linear programming
$\cdot$ mixed-integer nonlinear programming
$\cdot$ optimization solver
$\cdot$ branch-and-cut
$\cdot$ branch-and-price
$\cdot$ column generation
$\cdot$ parallelization
$\cdot$ mixed-integer semidefinite programming

\paragraph{\bf Mathematics Subject Classification} 90C05 $\cdot$ 90C10 $\cdot$ 90C11 $\cdot$ 90C30 $\cdot$ 90C90 $\cdot$ 65Y05

\newpage

\section{Introduction}
\label{sect:introduction}

The \scipopt comprises a set of complementary software packages designed to model and
solve a large variety of mathematical optimization problems:
the modeling language \zimpl~\cite{Koch2004},
the presolving library \papilo,
the linear programming solver
  \soplex~\cite{Wunderling1996},
the constraint integer programming solver
  \scip~\cite{Achterberg2009}, which can be used as a fast standalone global solver for
  mixed-integer linear and nonlinear programs and a flexible
  branch-cut-and-price framework,
the automatic decomposition solver \gcg~\cite{GamrathLuebbecke2010}, and
the \ug framework for solver parallelization~\cite{Shinano2018}.

All six tools can be downloaded in source code and are freely available
for members of noncommercial and academic institutions.
Development and bugfix branches of \scip, \soplex and \papilo are mirrored under
\url{https://github.com/orgs/scipopt} on a daily basis.
They are accompanied by several extensions for solving specific problem classes
such as the award-winning Steiner tree solver
\scipjack~\cite{Gamrath2017scipjack} and the mixed-integer semidefinite
programming (\MISDP) solver \scipsdp~\cite{GallyPfetschUlbrich2018}.
This paper discusses the capacity of \scip as a software and research tool and presents the
evolving possibilities for working with the \scipoptv, both as a black-box
toolbox and as a framework with possibilties of interaction and extension.

\paragraph{Background}

\scip is a branch-cut-and-price framework for
solving different types of optimization problems, most importantly,
\emph{mixed-integer linear programs} (\MILPs) and
\emph{mixed-integer nonlinear programs} (\MINLPs).
\MINLPs are optimization problems of the form
\begin{equation}
  \begin{aligned}
    \min \quad& \linobj^\T x \\
    \text{s.t.} \quad& \linmatrix x \geq \rhs, \\
    &\low{g}_k \leq g_k(x) \leq \upp{g}_k && \fa k \in \consindex, \\
    &\low{x}_{i} \leq x_{i} \leq \upp{x}_{i} && \fa i \in \varindex, \\
    &x_{i} \in \Z && \fa i \in \intvarindex,
  \end{aligned}
  \label{eq:generalmip}
\end{equation}
defined by $c \in \R^n$, $A \in\R^{m^{(\ell)}\times n}$, $ \rhs\in \R^{m^{(\ell)}}$,
$\low{g}$, $\upp{g}\in\Rinf^{m^{(n)}}$, $g : \R^{n} \rightarrow \R^{m^{(n)}}$,
$\low{x}$, $\upp{x} \in
\Rinf^{n}$, the index set of integer variables $\mathcal{I} \subseteq \mathcal{N} \defi \{1, \ldots, n\}$ and the index set of nonlinear constraints $\consindex \defi \{1,\ldots,m^{n}\}$.  We assume that $\nonlincons$ is specified in
algebraic form using basic expressions that are known to \scip.  The usage of $\Rinf \defi \R \cup
\{-\infty,\infty\}$ allows for variables that are free or bounded only in
one direction (we assume that no variable is fixed to~$\pm \infty$).
In the absence of nonlinear constraints $\low{g} \leq g(x) \leq \upp{g}$, the problem becomes an \MILP.

\scip is not restricted to solving MI(N)LPs, but is a framework for solving \emph{constraint integer programs} (\CIPs),
a generalization of the former two problem classes.
The introduction of \CIPs was motivated by the modeling flexibility of
constraint programming and the algorithmic requirements of integrating it with
efficient solution techniques available for \MILPs. Later on, this framework
allowed for the integration of \MINLPs.
Roughly speaking, \CIPs are finite-dimensional optimization problems with arbitrary constraints and a linear objective function
that satisfy the following property: if all integer variables are fixed, the remaining subproblem must form a linear or nonlinear program.

The core of \scip coordinates a central branch-cut-and-price algorithm that is augmented by a
collection of plugins.
The methods for processing constraints of a given type are implemented in 
\emph{constraint handler} plugins.
The default plugins included in the \scipopt provide tools to solve MI(N)LPs as
well as some problems from constraint programming, satisfiability testing and
pseudo-Boolean optimization.
In this way, advanced methods like primal heuristics, branching rules, and cutting plane separators can be integrated using a pre-defined interface.
\scip comes with many such plugins that enhance MI(N)LP performance, and new
plugins can be created by users.
This design and solving process is described in more detail by
Achterberg~\cite{Achterberg2007a}.

The core solving engine also includes \presollib, which provides an additional
presolving procedure that is called by \scip, and the linear programming (\LP)
solver \soplex which is used by default for solving the \LP relaxations within the
branch-cut-and-price algorithm.
Interfaces to several external \LP solvers exist, and new ones can be added by users.

The flexibility of this framework and its design, which is centered around the capacity for
extension and customization, are aimed at making \scip a versatile tool to be used by
optimization researchers and practitioners.
The possibility to modify the solving process by including own solver components enables users to
test their techniques within a general-purpose branch-cut-and-price
framework.

The extensions of \scip that are included in the \scipopt showcase the use of \scip as a basis for
the users' own projects.
\gcg extends \scip to automatically detect problem structure
and generically apply decomposition algorithms based on the Dantzig-Wolfe or the
Benders' decomposition scheme.
\scipsdp allows to solve mixed-integer semidefinite programs, and \scipjack is a solver for Steiner
tree problems.
Finally, the default instantiations of the \ug framework use \scip as a base
solver in order to perform branch-and-bound in parallel computing
environments.

\paragraph{Examples of Works Using \scip}

A number of works independent of the authors of this paper have employed \scip as a research tool.
Examples of such works include papers on new symmetry handling
algorithms~\cite{dias2021exploiting}, branching rules~\cite{anderson2021further} and integration of
machine learning with branch-and-bound based \MILP solvers~\cite{prouvost2021ecole}.
Further application-specific algorithms have been developed based on \scip, for example,
specialized algorithms for solving electric vehicle routing~\cite{ceselli2021branch} and network
path selection~\cite{casazza2021optimization} problems.
Many articles employ \scip as an \MINLP solver for problems such as hyperplanes
location~\cite{blanco2021multisource}, airport capacity extension, fleet investment, and optimal
aircraft scheduling~\cite{coniglio2021airport}, cryptanalysis problems~\cite{florez2021internal},
Wasserstein distance problems~\cite{ccelik2021wasserstein}, and chance-constrained nonlinear
programs~\cite{kannan2021stochastic}.

\paragraph{Structure of the Paper}

The paper is organised as follows.
A performance evaluation of \scipv and a comparison of its performance to that of \scip 7.0 is
carried out in Section~\ref{sect:performance}.
The core solving engine is discussed in Section~\ref{sect:core}.
The interfaces and modeling languages are presented in
Section~\ref{sect:languagesinterfaces}.
\scip extensions that are included in the \scipopt are discussed in Section~\ref{sect:extensions}, 
and Section~\ref{sect:finalRemarks} concludes the paper.

For a more detailed description of the new features introduced in \scipoptv, and for the technical
details, we refer the reader to the \scipoptv release report~\cite{BestuzhevaEtal2021ZR}.

\section{Performance of \scipv for MILP and MINLP}
\label{sect:performance}

In this section, we present computational experiments conducted by running \scip without parameter tuning or
algorithmic variations to assess the performance changes since the 7.0 release.
The indicators of interest
are
the number of solved instances, the shifted geometric mean
of the number of branch-and-bound nodes (shift 100~nodes), and the shifted geometric mean
of the solving time (shift 1~second).

\subsection{Experimental Setup}

We use the \scipopt~7.0 as the baseline, including \soplex~5.0 and \papilo~1.0,
and compare it with the \scipoptv including \soplexv and \papilov.
Both were compiled using GCC~7.5, use \solver{Ipopt}~3.12.13 as NLP subsolver built with the \solver{MUMPS}~4.10.0
numerical linear algebra solver, \solver{CppAD}~20180000.0 as algorithmic differentiation library,
and \solver{bliss}~0.73 for detecting symmetry.
The time limit was set to 7200 seconds in all cases.

The \MILP instances are selected from the \miplib 2003, 2010, and 2017~\cite{miplib2017} as well as the
\coral\cite{linderoth2005noncommercial} 
instance sets and include all
instances solved by \scip~7.0 with at least one of five random
seeds or solved by \scipv with at least one of five random seeds;
this amounts to 347~instances.
The \MINLP instances are similarly selected from the \minlplibtwo\footnote{\url{https://www.minlplib.org}}
with newly solvable instances added to the ones solved by \scip~7.0 for a total of 113 instances.

All performance tests were run on identical machines with Intel Xeon CPUs E5-2690 v4 @ 2.60GHz
and 128GB in RAM. A single run was carried out on each machine in a single-threaded mode.
Each optimization problem was solved with SCIP using five different seeds for random number generators.
This results in a testset of 565 \MINLPs and 1735 \MILPs.
Instances for which the solver reported numerically inconsistent results are excluded from the presented results.

\subsection{\MILP Performance}

Results of the performance runs on \MILP instances are presented in Table~\ref{tbl:rubberband_table_mip}.
The ``affected'' subset contains instances for which the two solver versions show different numbers of dual simplex iterations.
Instances in the subsets $[t,\texttt{tilim}]$ were solved by
at least one solver version within the time limit and took least $t$~seconds to solve with at least one version.
``\alloptimal'' and ``\difftimeouts'' are the subsets of instances that can
be solved by both versions and by exactly one version, respectively. ``relative'' shows the ratio of the shifted geometric mean
between the two versions.

The changes introduced with
\scip 8.0 improved the performance on \MILPs both in terms of number of solved instances
and time.
The improvement is more limited on `\alloptimal' instances that were solved by both solvers, for which the relative improvement is only of $12\,\%$.
This indicates that the overall speedup is more due to newly solved instances than to improvement on instances that were already solved by \scip 7.0.

\begin{table}
\caption{Performance comparison for MILP instances}
\label{tbl:rubberband_table_mip}
\scriptsize

\begin{tabular*}{\textwidth}{@{}l@{\;\;\extracolsep{\fill}}rrrrrrrrr@{}}
\toprule
&           & \multicolumn{3}{c}{\scip~8.0+\soplex~6.0} & \multicolumn{3}{c}{\scip~7.0+\soplex~5.0} & \multicolumn{2}{c}{relative} \\
\cmidrule{3-5} \cmidrule{6-8} \cmidrule{9-10}
Subset                & instances &                                   solved &       time &        nodes &                                   solved &       time &        nodes &       time &        nodes \\
\midrule
\cleaninst            &      1708 &                                     1478 &      231.3 &         3311 &                                     1445 &      271.3 &         4107 &        1.17 &          1.24 \\
\affected             &      1475 &                                     1424 &      173.8 &         2843 &                                     1391 &      209.7 &         3611 &        1.21 &          1.27 \\
\cmidrule{1-10}
\bracket{0}{tilim}    &      1529 &                                     1478 &      154.4 &         2512 &                                     1445 &      184.6 &         3167 &        1.20 &          1.26 \\
\bracket{1}{tilim}    &      1470 &                                     1419 &      185.9 &         2870 &                                     1386 &      223.8 &         3647 &        1.20 &          1.27 \\
\bracket{10}{tilim}   &      1361 &                                     1310 &      248.1 &         3612 &                                     1277 &      303.1 &         4661 &        1.22 &          1.29 \\
\bracket{100}{tilim}  &      1000 &                                      949 &      537.1 &         7270 &                                      916 &      702.6 &        10262 &        1.31 &          1.41 \\
\bracket{1000}{tilim} &       437 &                                      386 &     1566.2 &        17973 &                                      353 &     2383.1 &        31707 &        1.52 &          1.76 \\
\difftimeouts         &       135 &                                       84 &     2072.7 &        19597 &                                       51 &     5062.1 &        69354 &        2.44 &          3.54 \\
\alloptimal           &      1394 &                                     1394 &      119.9 &         2048 &                                     1394 &      133.8 &         2330 &        1.12 &          1.14 \\
\bottomrule
\end{tabular*}
\end{table}


\subsection{\MINLP Performance}
\label{sect:perfminlp}

With the major revision of the handling of nonlinear constraints,
the performance of \scip on \MINLPs has changed considerably compared to \scip~7.0.
The results are summarized in Table~\ref{tbl:rubberband_table_minlp}.
More instances are solved by \scipv than by \scip~7.0,
and \scipv solves the instances for each of these subsets with a shorter shifted geometric mean time.
On the 386 instances solved by both versions, \scipv requires fewer nodes and less time.
The number of instances solved by only one of the two versions (\difftimeouts) is much higher than
reported in previous release reports with similar experiments, with 68 instances solved only by \scipv
and 49 instances solved only by \scip~7.0.
A performance evaluation that focuses only on the changes in handling nonlinear constraints is given in Section~\ref{sect:perfconsexpr}.

\begin{table}
\caption{Performance comparison for MINLP}
\label{tbl:rubberband_table_minlp}
\scriptsize

\begin{tabular*}{\textwidth}{@{}l@{\;\;\extracolsep{\fill}}rrrrrrrrr@{}}
\toprule
&           & \multicolumn{3}{c}{\scip~8.0+\soplex~6.0} & \multicolumn{3}{c}{\scip~7.0+\soplex~5.0} & \multicolumn{2}{c}{relative} \\
\cmidrule{3-5} \cmidrule{6-8} \cmidrule{9-10}
Subset                & instances &                                   solved &       time &        nodes &                                   solved &       time &        nodes &       time &        nodes \\
\midrule
\cleaninst            &       558 &                                      454 &       39.1 &         2427 &                                      435 &       45.7 &         1845 &        1.17 &          0.76 \\
\affected             &       487 &                                      438 &       23.5 &         1748 &                                      419 &       28.4 &         1456 &        1.21 &          0.83 \\
\cmidrule{1-10}
\bracket{0}{tilim}    &       503 &                                      454 &       21.7 &         1585 &                                      435 &       25.9 &         1326 &        1.19 &          0.84 \\
\bracket{1}{tilim}    &       375 &                                      326 &       56.1 &         3994 &                                      307 &       71.0 &         3113 &        1.27 &          0.78 \\
\bracket{10}{tilim}   &       293 &                                      244 &      121.6 &         7450 &                                      225 &      169.3 &         5393 &        1.39 &          0.72 \\
\bracket{100}{tilim}  &       195 &                                      146 &      307.6 &        14204 &                                      127 &      433.9 &         6696 &        1.41 &          0.47 \\
\bracket{1000}{tilim} &       153 &                                      104 &      466.9 &        23425 &                                       85 &      565.3 &         8382 &        1.21 &          0.36 \\
\difftimeouts         &       117 &                                       68 &      451.4 &        29142 &                                       49 &      461.8 &         6275 &        1.02 &          0.22 \\
\alloptimal           &       386 &                                      386 &        8.2 &          609 &                                      386 &       10.4 &          806 &        1.27 &          1.32 \\
\bottomrule
\end{tabular*}
\end{table}

\section{The Core Solving Engine}
\label{sect:core}

This section presents the core solving engine, which includes the \CIP solver \scip, the \MILP
presolving library \papilo, and the \LP solver \soplex.
It discusses \scip's \MINLP framework in Section~\ref{sect:minlp}, which was completely reworked
in the 8.0 release, and demonstrates the possibilities for implementing user's own methods
using the examples of two areas that saw improvement with the 8.0 release, namely symmetry
handling and primal heuristics in Sections~\ref{sec:symmetry}
and~\ref{subsect:heuristics}.

The full list of new features introduced in \scipv is the following: a new framework for handling nonlinear constraints,
symmetry handling on general variables and improved orbitope detection,
a new separator for mixing cuts,
improvements to decomposition-based heuristics,
the option to apply the mixed integer rounding procedure when  generating optimality cuts in the Benders' decomposition framework,
a new plugin type that enables users to include their own cut selection rules into SCIP, and
several technical improvements.

Further, the section provides an overview of the presolving library \papilo and the \LP solver
\soplex in Sections~\ref{subsect:papilo} and~\ref{subsect:soplex}, and presents the new dual
postsolving feature in \papilo, which allowed for it to be integrated into \soplex.

\subsection{SCIP's New MINLP Framework}\label{sect:minlp}

A new framework for handling nonlinear constraints was introduced with the \scipv release.
The main motivation for this change is twofold: First, it aims at increasing the reliability of the
solver and alleviating numerical issues that arose from problem reformulations.
Second, the new design of the nonlinear framework reduces the ambiguity of expression and structure
types by implementing different kinds of plugins for low-level expressions that define expressions,
and high-level structures that add functionality for particular, often overlapping structures.

The main components of the new framework are the following:
plugins representing expressions;
a reimplementation of the constraint handler for nonlinear constraints, \texttt{cons\_nonlinear};
nonlinear handler plugins that provide functionality for high-level structures;
a revision of the primal heuristic that solves NLP subproblems;
revised interfaces to NLP solvers; and
revised interface to an automatic differentiation library.
Moreover, \scipv contains cutting plane separators that work on nonlinear structures and interact with
\texttt{cons\_nonlinear}.

\subsubsection{New Expressions Framework}
\label{sect:expr}

Algebraic expressions are well-formed combinations of constants, variables, and various algebraic operations such as addition, multiplication, exponentiation, that are used to describe mathematical functions.
In SCIP, they are represented by a directed acyclic graph with nodes representing variables, constants, and operators and arcs indicating the flow of computation.

With \scipv, the expression system has been completely rewritten.
Proper SCIP plugins, referred to as \textit{expression handlers}, are now used to define all semantics of an operator.
These expression handlers support more callbacks than what was available for user-defined operators before.
Furthermore, much ambiguity and complexity is avoided by adding expression handlers for basic operations only.
High-level structures such as quadratic functions can still be recognized, but are no longer made explicit by a change in the expression type.

\subsubsection{New Handler for Nonlinear Constraints}
\label{sect:consnl}

For \scipv, the constraint handler for nonlinear constraints, \texttt{cons\_nonlinear}, has been rewritten and constraint handlers for quadratic, second-order cone, absolute power, and bivariate constraints have been removed.
Some functionalities of the removed constraint handlers have been reimplemented in other plugins.

An initial motivation for rewriting \texttt{cons\_nonlinear} was a numerical issue which was caused by explicit constraint reformulation in earlier versions.
Such a reformulation can lead to a difference in constraint violation estimation in the
original and reformulated problems and, in particular, to a solution being feasible for the
reformulated problem and infeasible for the original problem.
For example, this occurs in a problem where the constraint $\exp(\ln(1000)+1+x\,y) \leq z$ is
reformulated as $\exp(w) \leq z$, $\ln(1000)+1+x\,y = w.$
On the MINLPLib library, this issue occurred for 7\% of instances.

The purpose of the reformulation is to enable constructing a linear relaxation.
In this process, nonlinear functions are approximated by linear under- and overestimators.
Since the formulas that are used to compute these estimators are only available for ``simple'' functions, new variables and constraints were introduced to split more complex expressions into adequate form~\cite{SmithPantelides1999,VigerskeGleixner2017}.

A trivial attempt to solve the issue of solutions not being feasible in the original problem would have been to add a feasibility check before accepting a solution.
However, if a solution is not feasible, actions to resolve the violation of original constraints need to be taken, such as a separating hyperplane, a domain reduction, or a branching operation.
Since the connection from the original to the presolved problem was not preserved, it would not have been clear which operations on the presolved problem would help best to remedy the violation in the original problem.

Thus, the new constraint handler aims to preserve the original constraints by applying only transformations that, in most situations, do not relax the feasible space when taking tolerances into account.
The reformulations that were necessary for the construction of a linear relaxation are not applied explicitly anymore, but handled implicitly by annotating the expressions that define the nonlinear constraints.
Another advantage of this approach is a clear distinction between the variables that were present in the original problem and the variables added for the reformulation.
With this information, branching is avoided on variables of the latter type.
Finally, it is now possible to exploit overlapping structures in an expression simultaneously.

\subsubsection{Extended Formulations}
Consider problems of the form \eqref{eq:generalmip}, where the set of nonlinear constraints is non-empty, and some constraints may be nonconvex. \scip solves such problems to global optimality via a spatial branch-and-bound algorithm.
Important parts of the algorithm are presolving, domain propagation, linear relaxation, and branching.
For domain propagation and linear relaxation, extended formulations are used which
are obtained by introducing \emph{slack variables} and replacing sub-trees of the expressions that define nonlinear constraints by \emph{auxiliary variables}.

These extended formulations have the following form:
\begin{equation}
  \label{eq:minlp_extdp}
  \tag{$\text{MINLP}_\text{ext}$}
  \begin{aligned}
    \min\quad & \linobj^\T x, \\
    \mathrm{s.t.}\quad & h_i(x,w_{i+1},\ldots,w_{m}) = w_i, & i=1,\ldots,m, \\
    & \low{x} \leq x \leq \upp{x},
    ~ \low{w} \leq w \leq \upp{w},
    ~ x_\intvarindex \in \Z^{\intvarindex}.
  \end{aligned}
\end{equation}
Here, $w_1,\ldots,w_m$ are slack variables, and $h_i\defi g_i$ for $i=1,\ldots,m$.
For each function $h_i$, subexpressions $f$ may be replaced by new auxiliary variables $w_{i'}$, $i'>m$, and new constraints $h_{i'}(x) = w_{i'}$ with $h_{i'} \defi f$ are added.
For the latter, subexpressions may be replaced again.
The result is referred to by $h_i(x,w_{i+1},\ldots,w_m)$ for any $i=1,\ldots, m$.
That is, to simplify notation, $w_{i+1}$ is used instead of $w_{\max(i,m)+1}$.

\paragraph{Example of an Extended Formulation}
Consider constraint
$
  \log(x)^2 + 2\log(x)\,y+y^2 \leq 4.
$
SCIP may replace $\log(x)$ by an auxiliary variable $w_2$, since this results in a quadratic form $w_2^2+2w_2y+y^2$, which is both bivariate and convex, the former being well suited for domain propagation and the latter being beneficial for linearization.
Therefore, the following extended formulation may be constructed:
\begin{align*}
&h_1(x,y,w_2) \defi (w_2)^2 + 2w_2y + y^2  = w_1,
 \\
&h_2(x,y) \defi \log(x)  = w_2,
~w_1 \leq 4.
\end{align*}

\subsubsection{Structure Handling}
\label{sect:nlhdlr}

The construction of extended formulations is based on the information on what algorithms are available for analyzing expressions of a specific structure.
Following the spirit of the plugin-oriented design of SCIP, these algorithms are added as separate plugins, referred to as \emph{nonlinear handlers}.
Nonlinear handlers can detect structures in expressions and provide domain propagation and linear relaxation algorithms that act on these structures.
Unlike other plugins in SCIP, nonlinear handlers are managed by \texttt{cons\_nonlinear} and not the SCIP core.

Nonlinear handlers for the following expression types are available in SCIP:
quadratic expressions defined as sums where at least one term is either a product of two expressions or a square expression,
bilinear expressions, convex and concave expressions, quotient expressions of the form $(ay_1+b)/(cy_2+d)+e$, and
expressions defined in terms of semi-continuous variables.
The second-order cone (SOC) nonlinear handler provides separation for SOC constraints.
Finally, the default nonlinear handler ensures that there always exist domain propagation and linear
under/overestimation routines for an expression and employs callbacks of expression handlers to provide the necessary functionalities.

Additional structures can be recognized for generating cutting planes to
strengthen \LP relaxations.
Such structures are handled by separator plugins.
While separators are not restricted to nonlinear structures, the following separators were
introduced in \scipv that work on \MINLPs:
the Reformulation-Linearization technique
(RLT)~\cite{adams1986tight,adams1990linearization,adams1993mixed} separator adds RLT cuts for
bilinear products and can additionally reveal linearized products between binary and continuous
variables;
the principal minor separator works on a matrix $X = xx^\T$, where entries
$X_{ij}$ represent auxiliary variables corresponding to $x_ix_j$, and enforces that principle
$2 \times 2$ minors are PSD;
and the intersection cuts separator for rank-1 constraints (disabled by default) adds
cuts derived from the condition that any $2\times 2$ minor of $X$ has determinant 0.

\subsubsection{Performance Impact of Updates for Nonlinear Constraints}
\label{sect:perfconsexpr}

While Section~\ref{sect:perfminlp} compared the performance of SCIP 7.0 and SCIP 8.0, this section takes a closer look at the effect of replacing only the handling of nonlinear constraints in SCIP.
That is, here the following two versions of SCIP are compared:
\begin{description}
\item \textbf{classic}: the main development branch of SCIP as of 23.08.2021; nonlinear constraints handled as in SCIP 7.0;
\item \textbf{new}: as classic, but with the handling of nonlinear constraints replaced as detailed in this section and symmetry detection extended to handle nonlinear constraints (see Section~\ref{sec:symmetry}).
\end{description}

SCIP has been build with GCC 7.5.0 and uses \papilo~1.0.2, \bliss~0.73, \cplex~20.1.0.1 as LP solver, \ipopt~3.14.4, \cppad~20180000.0 and Intel MKL 2020.4.304 for linear algebra (LAPACK).
\ipopt uses the same LAPACK and HSL MA27 as linear solver.
All runs are carried out on machines with Intel Xeon CPUs E5-2660 v3 @ 2.60GHz and 128GB RAM in a single-threaded mode.
A time limit of one hour, a memory limit of 100000MB, an absolute gap tolerance of $10^{-6}$, and a relative gap tolerance of $10^{-4}$ are set.
All 1678 instances of \minlplibtwo (version 66559cbc from 2021-03-11) that can be handled by both versions are used.
Note that \minlplibtwo is not designed to be a benchmark set, since, for example, some models are overrepresented. 
For each instance, two additional runs were conducted where the order of variables and constraints were permuted.
Thus, in total 5034 jobs were run for each version.

\begin{table}
 \centering
 \caption{Comparison of performance of SCIP with classic versus new handling of nonlinear constraints on MINLPLib.}
 \label{tab:minlplib_status}
 \small
 \begin{tabular}{@{}lrlrrr@{}}
  \toprule
  Subset & instances & metric & classic & new & both \\ \midrule
  all & 5034 & solution infeasible & 481 & 49 & 20 \\
    & & failed & 143 & 70 & 18 \\
    & & solved & 2929 & 3131 & 2742 \\
    & & time limit & 1962 & 1833 & 1598 \\
    & & memory limit & 0 & 0 & 0 \\
  clean & 4839 & fastest & 3733 & 3637 & 2531 \\
    & & mean time & 75.9s & 70.3s \\
    & & mean nodes & 2543 & 2601 \\
  \bottomrule
 \end{tabular}
\end{table}

Table~\ref{tab:minlplib_status} summarizes the results.
A run is considered as failed if the reported primal or dual bound conflicts with best known bounds for the instance, the solver aborted prematurely due to a fatal error,
or the solver did not terminate at the time limit.
Runs where the final solution is not feasible are counted separately.
With the new version, for much fewer instances the final incumbent is not feasible for the original problem, that is, the issue discussed in Section~\ref{sect:consnl} has been resolved.
For the remaining 49 instances, typically small violations of linear constraints or variable bounds occur.
Furthermore, the reduction in ``failed'' instances by half shows that the new version is more robust regarding the computation of primal and dual bounds.
Finally, the new version solves about 400 additional instances in comparision to the classic one, but also no longer solves about 200 instances within the time limit.

Subset ``clean'' refers to all instances where both versions did not fail, i.e., either solved to optimality or stopped due to the time limit.
We count a version to be ``fastest'' on an instance if it is not more than 25\% slower than the other version.
Mean times were computed as explained in the beginning of Section~\ref{sect:performance}.
Due to the increase in the number of solved instances, a reduction in the mean time with the new version on subset ``clean'' can be observed, even though the new version is fastest on less instances than the classic one.

\begin{figure}
 \begin{minipage}{.49\textwidth}
 \centering
 \includegraphics{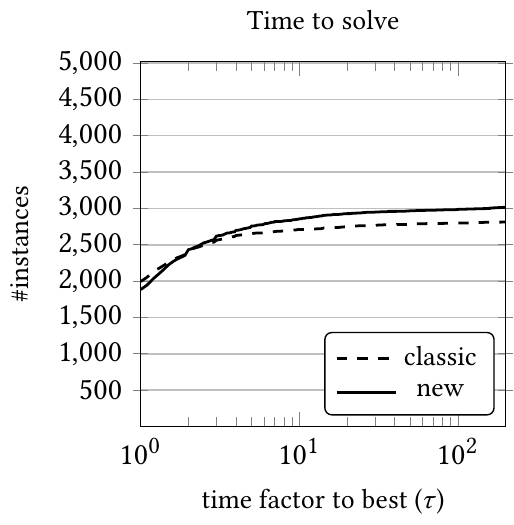}
 \end{minipage}
 \begin{minipage}{.49\textwidth}
 \centering
 \includegraphics{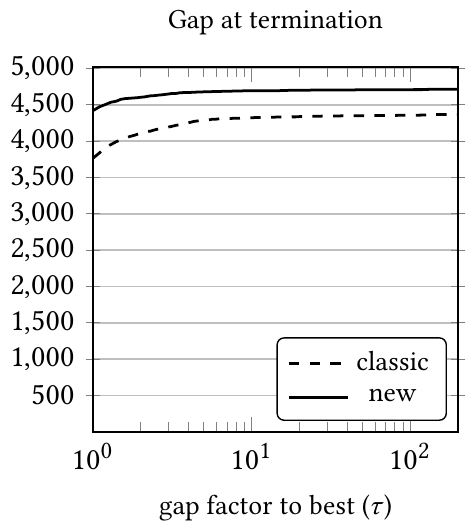}
 \end{minipage}

 \caption{Performance profiles with classic versus new handling of nonlinear constraints, showing the number of instances for which the corresponding version was at most $\tau$ times worse (regarding time (left) or gap at termination (right)) than the best of both versions.
 For the time plot, instances that were solved to optimality are considered.
 For the gap plot, instances that did not fail are considered.}
 \label{fig:minlplib_perfprofile}
\end{figure}

Figure~\ref{fig:minlplib_perfprofile} shows performance profiles that compare both versions w.r.t.\ the time to solve an instance and the gap at termination.
The time comparison visualizes what has been observed in Table~\ref{tab:minlplib_status}: the new version solves more instances, but can be slower.
The gap comparison shows that on instances that are not solved, often the new version produces a smaller optimality gap than the classic version.

\subsection{Improvements in Symmetry Handling}
\label{sec:symmetry}

Symmetries are known to have an adverse effect on the performance of
MI(N)LP solvers due to symmetric subproblems being treated
repeatedly without providing new information to the solver.
Since detecting all symmetries is \NP-hard~\cite{Margot2010}, \scip only
detects symmetries that keep the formulation invariant.

\scip's symmetry handling framework can be used both
as a black box and research tool.
In the black box approach, \scip automatically detects and handles symmetries.
If symmetries are known, users can tell \scip about them by
adding specialized constraints.
Customized code can include such constraints via API
functions, but also black box \scip can be informed about symmetries via
parsing them from files in \scip's \CIP format.
Moreover, \scip facilitates research on symmetries as it stores all
symmetry information centrally in the \texttt{symmetry} propagator and
provides implementations of basic symmetry
operations such as stabilizer computations.

For a permutation~$\perm$ of the variable index set~$\{1,\dots,n\}$ and a
vector~$x \in \R^n$, we define~$\perm(x) = (x_{\inv\perm(1)}, \dots,
x_{\inv\perm(n)})$.
We say that~$\perm$ is a \emph{symmetry} of~\eqref{eq:generalmip} if the
following holds: $x \in \R^n$
is feasible for~\eqref{eq:generalmip} if and only if~$\perm(x)$ is feasible, and~$c^\T x =
c^\T\perm(x)$.
The set of all symmetries forms a group~$\group$, the \emph{symmetry group}
of~\eqref{eq:generalmip}.
If~$\group$ is a product group~$\group = \group_1 \otimes \dots
\otimes \group_k$, the variables affected by one factor of~$\group$
are not affected by any other factor.
In this case, \scip can apply different symmetry
handling methods for each factor.
The sets of all variables affected by a single factor are called components.

\scip~7.0 was only able to handle symmetries of binary variables in
MILPs using two paradigms:
a constraint-based approach or the pure propagation-based approach \emph{orbital
  fixing}~\cite{Margot2002,Margot2003,OstrowskiEtAl2011}.
For a symmetry~$\perm$, the constraint-based approach enforces that the
variable vector~$x$ is lexicographically not smaller than~$\perm(x)$.
This is implemented via three different constraint
handler plugins.
For single permutations~$\perm$, the \emph{symresack} and \emph{orbisack}
constraint handlers use separation and propagation~\cite{HojnyPfetsch2019}
techniques for enforcing the lexicographic requirement, also
c.f.~\cite{KaibelLoos2011}.
Additionally, if an entire factor~$\group_i$ of~$\group$ has a special
structure,
the \emph{orbitope} constraint handler applies specialized
techniques~\cite{SCIP7}.

\scip~8.0 extends the symmetry handling framework.
First, it allows to detect symmetries in \MINLPs~\cite{Wegscheider2019}.
Second,
in \scip~8.0, symmetries of general variables can be handled by
inequalities derived from the Schreier-Sims table (SST
cuts)~\cite{Salvagnin2018,LibertiOstrowski2014}.
These inequalities are based on a list of
\emph{leaders}~${\ell_1},\dots,{\ell_k}$ together with
suitably defined \emph{orbits}~$O_1,\dots,O_k$, leading to inequalities
$x_{\ell_i} \geq x_j$, $j \in O_i$, $i \in \{1,\dots, k\}$.
Users have a high degree of flexibility to control the
selection of orbits
and can thus select the most promising
symmetry handling strategy.
Third, orbitope detection has been extended to also detect suborbitopes, i.e., parts of the
symmetry group that allow to apply orbitopes.
Since adding suborbitopes did not turn out to always be beneficial, \scip
adds suborbitopes according to a strategy
that can combine suborbitopes and SST cuts;
adding SST cuts can be controlled by a user via parameters.

Furthermore, \scip~8.0 contains
improvements of previously available methods.
First, if orbisack constraints interact with set packing or partitioning
constraints in a certain way, they are automatically upgraded to orbitopes.
This upgrade has been made more efficient.
Second, the running time of the separation routine of cover inequalities for
symresacks has been improved from quadratic to linear by using the observation
from~\cite{HojnyPfetsch2019} that
minimal cover inequalities for symresacks can be separated by merging
connected components of an auxiliary graph.
The new implementation exploits that its connected components
are either paths or cycles.
Finally, propagation routines of the symresack and orbisack constraint handler now find all
variable fixings that can be derived from local variable bound information.

\subsection{Primal Decomposition Heuristics}
\label{subsect:heuristics}

Most \MILPs have sparse constraint matrices for which a (bordered) block-diagonal form might be obtained by permuting the rows/columns of the matrix.
Identifying such a form allows for potentially rendering large-scale complex problems considerably more tractable.
Solution algorithms or heuristics can be designed exploiting the underlying structure and yielding smaller, easier problems.
In this sense, a so-called \emph{decomposition} identifies subsets of rows and columns that are only linked to each other via a set of linking rows and/or linking columns, but are otherwise independent.

A decomposition consisting of $k\in \N$ blocks is a partition
$
	\mathcal{D}\coloneqq (D^{\text{row}},D^{\text{col}})
	\text{ with }
	D^{\text{row}} \coloneqq (D^{\text{row}}_{1},\dots,D^{\text{row}}_{k},L^{\text{row}}), \ D^{\text{col}} \coloneqq (D^{\text{col}}_{1},\dots,D^{\text{col}}_{k},L^{\text{col}})
$
of the rows/columns of the constraint matrix $\linmatrix$ into $k+1$ pieces each,
whereby it holds for all $i\in D^{\text{row}}_{q_{1}}$, $j\in D^{\text{col}}_{q_{2}}$ that $a_{i,j} \neq 0$ implies $q_{1} = q_{2}$.
Rows $L^{\text{row}}$ and columns $ L^{\text{col}}$, which may be empty, are called \emph{linking rows and columns}, respectively.

In general, there is no unique way to decompose an \MILP, and different
decompositions might lead to different solver behaviors.
Users might be aware of decompositions and know which are most useful
for a specific problem.
Therefore, since version 7.0 it is possible to pass user decompositions to \scip~\cite{SCIP7}.
A decomposition structure can be created using the \scip API, assigning labels to variables and/or constraints, and calling automatic label computation procedures if necessary.
Alternatively, \scip also provides a file reader for decompositions in constraints.

In \scip~7.0, the Benders decomposition framework and the heuristic Graph Induced Neighborhood
Search were extended to exploit user-provided decompositions, and a first version of the heuristic
Penalty Alternating Direction Method (PADM)~\cite{Geissleretal2017,ScheweSchmidtWeninger2019} was
introduced.
\scipv comes with an improvement of PADM
and provides another decomposition heuristic Dynamic Partition Search (DPS)~\cite{BestuzhevaEtal2021ZR}.

\paragraph{Improvement of Penalty Alternating Direction Method}
PADM splits an \MINLP into several subproblems
according to a given decomposition $\mathcal{D}$ with linking variables only,
whereby the linking variables get copied and the differences are penalized.
Then the subproblems are alternatingly solved.
For faster convergence, the objective function of each subproblem
has been replaced by a penalty term, and this replacement can lead to arbitrarily bad solutions.
Therefore, PADM has been extended by the option
to improve a found solution by reintroducing the original objective function.

\paragraph{Dynamic Partition Search}
The new primal construction heuristic DPS requires a decomposition with linking constraints only.
The linking constraints and their sides are split
by introducing vectors $p_q \in \R^{L^{\text{row}}}$ for each block $q \in \{1,\dots,k\}$,
where $\R^{L^{\text{row}}}$ denotes the space of vectors with components indexed by $L^{\text{row}}$,
and requiring that the following holds:
\begin{equation}\label{eq:DecompPartition}
  \sum_{q=1}^{k} p_q = b_{[L^{\text{row}}]}.
\end{equation}

To obtain information on subproblem infeasibility and speed up the solving process,
the objective function is replaced by a weighted sum of slack variables $z_q \in \R^{L^{\text{row}}}_+$.
For penalty parameter $\lambda \in \R^{L^{\text{row}}}_{> 0}$,
each subproblem $q$ has the form
\begin{equation}
\label{eq:dpsblock}
\begin{aligned}
\min\quad & \lambda^\T z_q, \\
\mathrm{s.t.}\quad & A_{[D^{\text{row}}_{q},D^{\text{col}}_{q}]}\,x_{[D^{\text{col}}_{q}]} \geq \rhs_{[D^{\text{row}}_{q}]}, \\
& \low{x}_{i} \leq x_{i} \leq \upp{x}_{i} && \fa i \in \varindex \cap D^{\text{col}}_{q}, \\
& x_{i} \in \Z && \fa i \in \intvarindex \cap D^{\text{col}}_{q},\\
& A_{[L^{\text{row}},D^{\text{col}}_{q}]}\,x_{[D^{\text{col}}_{q}]} + z_q \geq p_{q},\\
& z_q \in \R^{L^{\text{row}}}_+.
\end{aligned}
\end{equation}

From~\eqref{eq:dpsblock}, it is apparent that the correct choice of $p_{q}$ plays a central role.
For this reason, we refer to $(p_q)_{q \in \{1,\dots,k\}}$ as a \emph{partition} of $b_{[L^{\text{row}}]}$.
The method starts with an initial partition fulfilling~\eqref{eq:DecompPartition}.
Then it is checked whether this partition will lead to a feasible solution
by solving $k$ independent subproblems~\eqref{eq:dpsblock}
with fixed $p_q$.
If the current partition does not correspond to a feasible solution, then the partition
gets updated, so that~\eqref{eq:DecompPartition} still holds.
These steps are repeated.
Similarly to PADM, it is possible to improve the found solution
by reoptimizing with the original objective function.

\subsection{\papilo}
\label{subsect:papilo}

The C++ library \papilo provides presolving routines for (MI)LP problems and was introduced with the \scipopt~7.0 \cite{SCIP7}.
\papilo can be integrated into \MILP solvers or used as a standalone presolver.
As a standalone presolver it provides presolving and postsolving routines.
Hence, it can be used to a) provide presolving for new solving methods and b) generate presolved instances so that different solvers can be benchmarked independently of their own presolvers.
Thus, the performance/behavior of the actual solver can be evaluated and compared more precisely.

\papilo's transaction-based design allows presolvers to run in parallel without requiring expensive copies of the problem and without special synchronizations.
Instead of applying results immediately, presolvers return their reductions to the core, where they are applied in a deterministic, sequential order.
Validity of every reduction to the modified problem is checked to avoid applying conflicting reductions.

Presolving deletes variables from the original problem by fixing, substituting, and aggregating variables. After solving the reduced problem, its solution does not contain any information on missing variables.
To restore the solution values of these variables and obtain a feasible solution of the original problem, corresponding data needs to be stored during the presolving process.
The process of recalculating the original solution from the reduced one is called postsolving or post-processing~\cite{AchterbergBixbyGuetal.2019}.
Until version 1.0.2, \papilo supported only postsolving primal solutions for LPs.
In the latest version, \papilo supports postsolving also for dual solutions, reduced costs, slack variables of the constraints, and the basic status of the variables and constraints for the majority of the LP presolvers.

\subsection{\soplex}
\label{subsect:soplex}

\soplex is a simplex-based \LP solver and an essential part of the optimization suite, since is the default \LP-solver for \scip.
In addition to all the essential features of a state-of-the-art \LP solver such as scaling, exploitation of sparsity, or presolving, \soplex also
supports an option for 80bit extended precision and an iterative refinement algorithm to produce high-precision solutions.
This enables \soplex to also compute exact rational solutions to \LPs, using either continued fraction approximations or a symbolic LU factorization.

The support of postsolving of dual LP solutions and basis information in \papilo makes it possible to integrate \papilo fully as a presolving library into \soplex.
In version~6.0 of \soplex, \papilo is available as an additional option for presolving.  The previous presolving implementation continues to be the default.

\section{Modeling Languages and Interfaces}
\label{sect:languagesinterfaces}

There are many interfaces to SCIP from different programming and modeling languages.
These interfaces allow users to programmatically call SCIP with an API close to the C one or leverage a higher-level syntax.

The AMPL interface has been rewritten and moved to the main SCIP library and executable.
With the \scipoptv, there exists a C wrapper for SoPlex, updated GAMS interfaces for SoPlex and
SCIP, a Julia package \texttt{SCIP.jl}, a basic Java interface JSCIPOpt, a new
Matlab interface for \scipv and \scipsdp based on the \texttt{OPTI Toolbox} by Jonathan Currie, and
the Python interface PySCIPOpt which can now also be installed as a Conda package.

The modeling language \zimpl allows for MI(N)LPs to be written and translated
into some file formats supported by \scip.
\zimpl 3.5.0 allows quadratic objective functions in addition to previously supported linear objective functions,
and can write suitable
instances as Quadratic Unconstrained Binary Optimization problems.

\section{Extensions}
\label{sect:extensions}

\subsection{The GCG Decomposition Solver}
\label{sect:gcg}

\scip allows implementing tailored decomposition-based algorithms.
Complementary to this, \gcg turns \scip into a generic decomposition-based \emph{solver\/}
for MILPs. While \gcg's focus is
on Dantzig-Wolfe reformulation (DWR) and Lagrangian decomposition,
Benders decomposition (BD) is also supported. The philosophy
behind \gcg is that decomposition-based algorithms
can be routinely applied to MILPs without the user's
interaction or even knowledge. To this end,
\gcg automatically detects a model structure that admits a
decomposition
and performs the corresponding
reformulation. 
This results in a master problem and one or several subproblems, which
are usually formulated as MILPs.
Based on the reformulation, the linear
relaxation in every node is solved by column generation (in the DWR
case) and Benders cut generation (in the BD case).  \gcg
features primal heuristics and separation of cutting planes, several
of which are adapted from \scip, but some are tailored to the
decomposition situation in which both an original and a reformulated
model are available.
As a research tool, \gcg can be used to quickly assess the potential
of a decomposition-based algorithm for any problem for which a compact
MILP formulation is available. 
This allows evaluating the performance of an algorithmic idea with a
single generic implementation, but across many different applications.
In what follows, we describe some enhancements in the
\gcg~3.5 release.

\subsubsection{Detection Loop Refactoring}
%

Decomposition-based algorithms rely on model structures,
cf.\ Section~\ref{subsect:heuristics}.
For automatic identification of such structures, \gcg features a modular detection loop.
\emph{Detectors}
iteratively assign roles like ``master'' or ``block'' to variables
and/or constraints.
This way, usually many different potential decompositions are
found.  We refer to the \scipopt~6.0 release report~\cite{SCIP6} for
a more detailed overview.  Detectors are implemented as plugins such that new ones can
be added conveniently. In every round, each detector works on existing
(but possibly empty) partial or complete decompositions. An empirically
very successful detection concept builds on the classification of
constraints and variables, which is performed prior to the actual
detection process, using \emph{classifiers}.

\subsubsection{Branching}
\label{sect:gcgstrongbranch}

In \gcg, two general branching rules are implemented (branching on
original variables~\cite{VilleneuveDesrosiersLuebbeckeSoumis:05} and
Vanderbeck's generic branching~\cite{vanderbeck11}) as well as one
rule that applies only to set partitioning master problems
(Ryan and Foster branching~\cite{RyanFoster81}).
While these rules differ significantly,
the general procedure has
two common stages: First, one determines the set of candidates we could
possibly branch on (called the branching rule).
Second, the \emph{branching candidate selection heuristic\/} 
selects one of the candidates.
\gcg previously contained pseudo cost, most fractional,
and random branching as selection heuristics for original variable
branching, and first-index branching for Ryan-Foster and Vanderbeck's
generic branching. In \gcg~\gcgversion, new strong branching-based selection heuristics are added~\cite{gaul2021branching}.


\subsubsection{Python Interface}
With \gcg~\gcgversion\
we introduce \textsc{PyGCGOpt} which extends \scip's existing Python
interface~\cite{MaherMiltenbergerPedrosoRehfeldtSchwarzSerrano2016}
for \gcg and is distributed independently from
the optimization suite\footnote{\url{https://github.com/scipopt/PyGCGOpt}}.
All existing functionality for MILP
modeling 
 is inherited from \textsc{PySCIPOpt}; therefore, any MILP
modeled in Python can be solved with \gcg without additional
effort.
The
interface supports specifying custom decompositions and exploring
automatically detected decompositions,
and plugins for {detectors} and {pricing solvers} can be implemented in Python.

\subsubsection{Visualization Suite}
\label{sec:gcg_visu_suite}

Visualizations of algorithmic behavior can yield understanding and
intuition for interesting parts of a solving process.
With \gcg~\gcgversion, we include a Python-based
\emph{visualization suite} that offers
visualization scripts to show processes and results related to
detection, branching, or pricing, among others. 
We highlight two features:


\begin{enumerate}
\item Reporting functionality:
  A
  \emph{decomposition report}
   offers an overview of all decompositions that \gcg found
  for a single run. For different runs, \gcg~3.5 offers two reports:
  A \emph{testset report} shows data and graphics
  for each single run of one selected test set. A \emph{comparison
    report} allows to compare two or more runs on the same test
  set.
\item Jupyter notebook: data produced for the reports can
  be 
  read, cleaned, filtered, and visualized interactively.
\end{enumerate}


\subsection{SCIP-SDP}
\label{sect:SCIP-SDP}

\scipsdp is an MISDP solver and a platform
for implementing methods for solving MISDPs. It was initiated by Sonja Mars and Lars Schewe~\cite{Mar13}, and then continued by Gally et
al.~\cite{GallyPfetschUlbrich2018} and Gally~\cite{Gal19}. New
results and methods mainly concerning presolving and
propagation are presented in \cite{MatP21}. \scipsdp features
interfaces to SDP-solvers DSDP, Mosek, and SDPA.

\scipsdp implements an SDP-based branch-and-bound method, which solves a
continuous SDP relaxation in each node. It incorporates
plugins such as primal heuristics, presolving and
propagation methods, and file readers.
There is also an option to solve LP relaxations in each node of the
branch-and-bound tree and generate
eigenvector cuts, see Sherali and Fraticelli~\cite{SheraliF2002}. This is sometimes faster than solving SDPs
in every node. These two options can also be run concurrently
if the parallel interface \texttt{TPI} of \scip is used. There also is a Matlab interface to
\scipsdp.

Moreover, \scipsdp can handle rank-1 constraints, that is, the
requirement that a matrix~$A$ has rank 1.
For such a constraint, quadratic constraints are added, modeling that
all $2 \times 2$-minors of~$A$ are zero~\cite{ChenAtamturkOren2017}.

Before we present some computational results, let us add some words of
caution. Although \scipsdp is numerically quite robust, accurately solving
SDPs is more demanding than solving LPs. This can lead to wrong results on
some instances\footnote{For instance, in seldom cases, the dual bound might
  exceed the value of a primal feasible solution.}, and the results often
depend on the tolerances.
Moreover,
the SDP-solvers use relative tolerances, while \scipsdp uses absolute
tolerances. Finally, for Mosek, we use a slightly tighter feasibility tolerance than in
\scipsdp.

Table~\ref{tab:scipsdp40vs32} shows a comparison between \scipsdp 3.2 and
4.0 on the same testset as used by Gally et
al.~\cite{GallyPfetschUlbrich2018}, which consists of 194 instances; the
changes between \scipsdp~4.0 and 3.2 are presented in more detail
in~\cite{BestuzhevaEtal2021ZR,MatP21}.
Reported are the number of optimally solved instances, as well as the
shifted geometric means of the number of processed nodes and the CPU time
in seconds. We use Mosek 9.2.40 for solving the continuous
SDP relaxations. The tests were performed on a Linux cluster with 3.5 GHz
Intel Xeon E5-1620 Quad-Core CPUs, having 32~GB main memory and 10~MB
cache. All computations were run single-threaded and with a timelimit of
one hour.

As can be seen from the results, \scipsdp~4.0 is considerably faster than
\scipsdp~3.2, but we recall that we have relaxed the tolerances.
Nevertheless, we conclude that \scipsdp has significantly
improved since the last version.

\begin{table}
  \caption{Performance comparison of \scipsdp 4.0 vs. \scipsdp 3.2}
  \label{tab:scipsdp40vs32}
  \small
  \begin{tabular*}{\textwidth}{@{}l@{\;\;\extracolsep{\fill}}ccc@{}}
    \toprule
    &  \# opt & \# nodes & time [s] \\
    \midrule
    \scipsdp 3.2 & 185 & 617.3 & 42.9 \\
    \scipsdp 4.0 & 187 & 497.3 & 26.6 \\
    \bottomrule
  \end{tabular*}
\end{table}

\subsection{\scipjack: Solving Steiner Tree and Related Problems}
\label{sect:SCIP-Jack}

Given an undirected, connected graph, edge costs and a set of \textit{terminal} vertices, the \textit{Steiner tree problem in graphs} (SPG)
asks for a tree of minimum weight that covers all terminals.
The SPG is a fundamental $\mathcal{NP}$-hard problem~\cite{Karp72} and one of the most studied problems in combinatorial optimization.

\scipjack, an exact SPG-solver, is built on the  branch-and-cut framework provided by \scip and makes extensive use of its plugin-based design.
At the heart of the implementation is a constraint handler that separates violated constraints, most importantly the so-called \emph{directed Steiner cuts}, which are separated by a specialized maximum-flow algorithm~\cite{Rehfeldt2021}.
The implementation includes a variety of additional \scip plugins, such as heuristics, propagators, branching rules and relaxators.
Finally, the use of \scip provides significant flexibility in the model to be solved, for example it is easily possible to add additional constraints. In this way, \scipjack can solve not only the SPG, but also 14 related problems.

The \scipoptv contains the new \scipjack~\scipjackversion\footnote{See also \url{https://scipjack.zib.de}.}, which comes with major improvements on most problem classes it can handle 
and outperforms the SPG solver by Polzin and Vahdati~\cite{PolzinThesis, VahdatiThesis}, which had remained unchallenged for almost 20 years, on almost all nontrivial benchmark testsets from the literature~\cite{RehfeldtKoch20Conflicts}.

\begin{figure}
\includegraphics[scale=0.8]{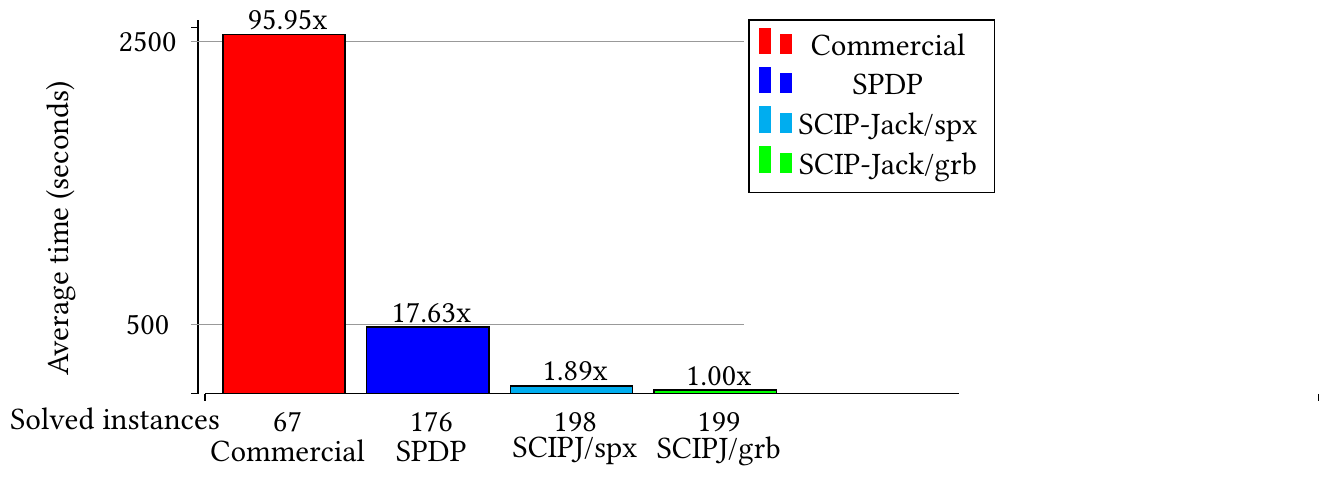}
 \caption{Computational results on the 200 benchmark instances of Tracks A and B of the PACE Challenge 2018.}\label{fig:stp:bra}
 \end{figure}

Figure~\ref{fig:stp:bra} provides computational results on the instances from Tracks A and B of the PACE Challenge 2018~\cite{BonnetS18}.
We use \gurobi 9.5 (\emph{Commercial}), the best other solver from the PACE Challenge (\emph{SPDP}~\cite{iwata2019}), and \scipjack with \soplex (\emph{SCIPJ/spx}) and
  \gurobi 9.5 (\emph{SCIPJ/grb}) as \LP solvers. A timelimit of one hour was set. Average times are given as arithmetic means with time-outs counted as one hour each.
The results were obtained on Intel Xeon CPUs E3-1245 @ 3.40~GHz with 32 GB RAM. It can be seen that SCIPJ/grb is roughly 17 times faster than SPDP, and 96 times faster than \gurobi.
For larger instances of the PACE 2018 benchmark, one commonly observes a run time difference of more than six orders of magnitude between \scipjack and commercial \MILP solvers.

Considerable problem-specific improvements have been made for the prize-collecting Steiner tree problem (STP) and (to a lesser extent) for the maximum-weight connect subgraph problem~\cite{RehfeldtKoch2021,RehfeldtFranzKoch2020}.
\scipjack \scipjackversion~can solve many previously unsolved benchmark instances from both problem classes to optimality---the largest of these instances have up to 10 million edges. 
Large improvements are observed for the Euclidean STP:
\scipjack \scipjackversion~is able to solve 19 Euclidean STPs with up to 100\,000 terminals to optimality for the first time~\cite{RehfeldtKoch20Conflicts}. Notably, the state-of-the-art Euclidean STP solver GeoSteiner 5.1~\cite{juhl2018geosteiner} could not solve any of these instances, even after one week of computation. In contrast, \scipjack \scipjackversion~solves all of them within 12 minutes, some even within two minutes.

\subsection{The UG Framework}
\label{sect:ug}

UG is a generic framework for parallelizing solvers in a distributed or shared memory computing environment. 
It was designed to parallelize state-of-the-art branch-and-bound solvers externally in order to exploit their powerful performance.
We have developed parallel solvers for SCIP~\cite{Shinano2012, Shinano18fiber, Shinano-ParaSCIP}, CPLEX (not developed anymore), FICO Xpress~\cite{Shinano-ParaXpress}, 
PIPS-SBB~\cite{Munguia2016,Munguia2019}, Concorde\footnote{\url{https://www.math.uwaterloo.ca/tsp/concorde.html}}, and QapNB~\cite{FujiiItoKimetal.2021}. 
Customized SCIP-based solvers such as SCIP-SDP and SCIP-Jack can be parallelized with minimal effort~\cite{Shinano2019-1}. 
The parallel version of SCIP-Jack solved several previously unsolved instances from SteinLib~\cite{Koch2001SteinLibAU} by using up to 43,000 cores~\cite{Shinano2019-2}.

In addition to the parallelization of these branch-and-bound base solvers, UG was used to develop MAP-SVP~\cite{MAP-SVP}, which is a solver for the Shortest Vector Problem (SVP), whose algorithm does not rely on branch-and-bound. 
For these applications, UG had to be adapted and modified.
Especially, the success of MAP-SVP, which updated several records of the SVP challenge\footnote{\url{http://latticechallenge.org/svp-challenge}}, motivated us to develop \emph{generalized UG}, in which all solvers developed so far can be handled by a single unified framework.
This has enabled UG 1.0 to serve as the basis for the parallel frameworks CMAP-LAP (Configurable Massively Parallel solver framework for LAttice Problems)~\cite{TateiwaShinanoYamamuraetal.2021} and CMAP-DeepBKZ~\cite{TateiwaShinanoYasudaetal.2021}.

\section{Final Remarks}
\label{sect:finalRemarks}

We discussed the functionality that the \scipopt offers optimization researchers,
and highlighted performance improvements and new functionality that was introduced in the \scipoptv.
The performance comparison of \scip 7.0 and \scipv showed a 17\,\% speed-up on both the \MILP and \MINLP testsets.
This was followed by a discussion of some aspects of the core solving engine of the \scipopt.
The new framework for handling nonlinear constraints was presented, which offers increased
reliability as well as improved handling of different types of nonlinearities that reduces type
ambiguity and extends support for implementing the handling of user-defined nonlinearities.
The use of \scip's flexible plugin-based structure for extending the solver with user methods
was demonstrated on the examples of new symmetry handling methods and primal decomposition
heuristics.
The framework that \scip provides for working on these methods was explained and the relevant
plugin types and other customization-enabling features were discussed, followed by the presentation
of new methods added in \scipv.

Further, we presented extensions built around \scip.
The semidefinite programming solver \scipsdp and the Steiner tree problem solver \scipjack provide
users of the \scipopt the functionality for solving more problem classes, the decomposition solver
\gcg offers a different solving approach, and the solver parallelization framework \ug enables the
use of branch-and-bound solvers, and in particular \scip, in parallel computing environments.
Moreover, these components of the \scipopt demonstrate how \scip's features can be leveraged in
creating new research projects which can extend beyond \scip's standard focus and approach.

\subsection*{Acknowledgements}

The authors want to thank all previous developers and contributors to the \scipopt and
all users that reported bugs and often also helped reproducing and fixing the bugs.
In particular, thanks go to Suresh Bolusani, Didier Ch\'{e}telat, Gregor Hendel, Gioni Mexi, Matthias Miltenberger,
Andreas Schmitt, Robert Schwarz, Helena V\"{o}lker, Matthias Walter, and Antoine Prouvoust and the Ecole team.
The Matlab interface was set up with the big help of Nicolai Simon.

\subsection*{Author Affiliations}

\hypersetup{urlcolor=black}
\newcommand{\myorcid}[1]{ORCID: \href{https://orcid.org/#1}{#1}}
\newcommand{\myemail}[1]{E-mail: \href{#1}{#1}}
\newcommand{\myaffil}[2]{{\noindent #1}\\{#2}\bigskip}

\small

\noindent KSENIA BESTUZHEVA, Zuse Institute Berlin, Department AIS$^2$T, Takustr.~7, 14195~Berlin, Germany; e-mail: bestuzheva@zib.de; \myorcid{0000-0002-7018-7099}

\noindent MATHIEU BESAN\c{C}ON, ZIB-AIS$^2$T; e-mail: besancon@zib.de; \myorcid{0000-0002-6284-3033}

\noindent WEI-KUN CHEN, School of Mathematics and Statistics, Beijing Institute of Technology, Beijing 100081, China; e-mail: chenweikun@bit.edu.cn; \myorcid{0000-0003-4147-1346} 

\noindent ANTONIA CHMIELA, ZIB-AIS$^2$T; e-mail: chmiela@zib.de; \myorcid{0000-0002-4809-2958} 

\noindent TIM DONKIEWICZ, RWTH Aachen University, Lehrstuhl f\"ur Operations Research, Kackertstr.~7, 52072~Aachen, Germany; e-mail: tim.donkiewicz@rwth-aachen.de; \myorcid{0000-0002-5721-3563} 

\noindent JASPER VAN DOORNMALEN, Eindhoven University of Technology, 
Department of Mathematics and Computer Science; e-mail: 
m.j.v.doornmalen@tue.nl; \myorcid{0000-0002-2494-0705} 

\noindent LEON EIFLER, ZIB-AIS$^2$T; e-mail: eifler@zib.de; \myorcid{0000-0003-0245-9344} 

\noindent OLIVER GAUL, RWTH Aachen University, Lehrstuhl f\"ur Operations Research, Kackertstr.~7, 52072~Aachen, Germany; e-mail: oliver.gaul@rwth-aachen.de; \myorcid{0000-0002-2131-1911} 

\noindent GERALD GAMRATH, ZIB-AIS$^2$T and I$^2$DAMO GmbH, Englerallee 19, 14195~Berlin, Germany; e-mail: gamrath@zib.de; \myorcid{0000-0001-6141-5937} 

\noindent AMBROS GLEIXNER, ZIB-AIS$^2$T and HTW Berlin; e-mail: gleixner@zib.de; \myorcid{0000-0003-0391-5903} 

\noindent LEONA GOTTWALD, ZIB-AIS$^2$T; e-mail: gottwald@zib.de; \myorcid{0000-0002-8894-5011} 

\noindent CHRISTOPH GRACZYK, ZIB-AIS$^2$T; e-mail: graczyk@zib.de; \myorcid{0000-0001-8990-9912} 

\noindent KATRIN HALBIG, Friedrich-Alexander Universit{\"a}t Erlangen-N{\"u}rnberg, Department of Data Science, Cauerstr.~11, 91058~Erlangen, Germany; e-mail: katrin.halbig@fau.de; \myorcid{0000-0002-8730-3447} 

\noindent ALEXANDER HOEN, ZIB-AIS$^2$T; e-mail: hoen@zib.de; \myorcid{0000-0003-1065-1651} 

\noindent CHRISTOPHER HOJNY, Eindhoven University of Technology, 
Department of Mathematics and Computer Science; e-mail: c.hojny@tue.nl; 
\myorcid{0000-0002-5324-8996} 

\noindent ROLF VAN DER HULST, University of Twente, Department of Discrete Mathematics and Mathematical Programming, P.O. Box 217, 7500 AE Enschede, The Netherlands; e-mail: r.p.vanderhulst@utwente.nl;
\myorcid{0000-0002-5941-3016}

\noindent THORSTEN KOCH, Technische Universit\"at Berlin, Chair of Software and Algorithms for Discrete Optimization, Stra\ss{}e des 17. Juni 135, 10623 Berlin, Germany, and ZIB, Department A$^2$IM, Takustr. 7, 14195~Berlin, Germany; e-mail: koch@zib.de; \myorcid{0000-0002-1967-0077} 

\noindent MARCO L\"UBBECKE, RWTH Aachen-Lehrstuhl f\"ur OR; e-mail: marco.luebbecke@rwth-aachen.de; \myorcid{0000-0002-2635-0522} 

\noindent STEPHEN J.~MAHER, University of Exeter, College of Engineering, Mathematics and Physical Sciences, Exeter, United Kingdom; e-mail: s.j.maher@exeter.ac.uk, \myorcid{0000-0003-3773-6882} 

\noindent FREDERIC MATTER, Technische Universit{\"a}t Darmstadt, Fachbereich Mathematik, Dolivostr.~15, 64293~Darmstadt, Germany; e-mail: matter@mathematik.tu-darmstadt.de; \myorcid{0000-0002-0499-1820} 

\noindent ERIK M\"UHMER, RWTH Aachen-Lehrstuhl f\"ur OR; e-mail: erik.muehmer@rwth-aachen.de; \myorcid{0000-0003-1114-3800} 

\noindent BENJAMIN M{\"U}LLER, ZIB-AIS$^2$T; e-mail: benjamin.mueller@zib.de; \myorcid{0000-0002-4463-2873} 

\noindent MARC E.~PFETSCH, TU Darmstadt, Fachbereich Mathematik; e-mail: pfetsch@mathematik.tu-darmstadt.de; \myorcid{0000-0002-0947-7193} 

\noindent DANIEL REHFELDT, ZIB-A$^2$IM; e-mail: rehfeldt@zib.de; \myorcid{0000-0002-2877-074X} 

\noindent STEFFAN SCHLEIN, RWTH Aachen-Lehrstuhl f\"ur OR; e-mail: steffan.schlein@rwth-aachen.de; \myorcid{0009-0003-8322-7226}

\noindent FRANZISKA SCHL{\"O}SSER, ZIB-AIS$^2$T; e-mail: schloesser@zib.de; \myorcid{0000-0002-3716-5031}

\noindent FELIPE SERRANO, ZIB-AIS$^2$T; e-mail: serrano@zib.de; \myorcid{0000-0002-7892-3951} 

\noindent YUJI SHINANO, ZIB-A$^2$IM; e-mail: shinano@zib.de; \myorcid{0000-0002-2902-882X} 

\noindent BORO SOFRANAC, ZIB-AIS$^2$T and TU Berlin; e-mail: sofranac@zib.de; \myorcid{0000-0003-2252-9469} 

\noindent MARK TURNER, ZIB-A$^2$IM and Chair of Software and Algorithms for Discrete Optimization, Institute of Mathematics, TU Berlin; e-mail: turner@zib.de; \myorcid{0000-0001-7270-1496} 

\noindent STEFAN VIGERSKE, GAMS Software GmbH, c/o ZIB-AIS$^2$T; e-mail: svigerske@gams.com; \myorcid{0009-0001-2262-0601}

\noindent FABIAN WEGSCHEIDER, ZIB-AIS$^2$T; e-mail: wegscheider@zib.de; \myorcid{0009-0000-8100-6751}

\noindent PHILIPP WELLNER; e-mail: p.we@fu-berlin.de; \myorcid{0009-0001-1109-3877}

\noindent DIETER WENINGER, Friedrich-Alexander Universität Erlangen-Nürnberg, Department of Data Science, Cauerstr.~11, 91058~Erlangen, Germany; e-mail: dieter.weninger@fau.de; \myorcid{0000-0002-1333-8591} 

\noindent JAKOB WITZIG, ZIB-AIS$^2$T; e-mail: witzig@zib.de; \myorcid{0000-0003-2698-0767}

\setlength{\bibsep}{0.25ex plus 0.3ex}
\bibliographystyle{abbrvnat}
\begin{small}
\bibliography{scipopt}
\end{small}

\end{document}